\documentclass[a4paper,12pt]{article}
\usepackage{amsmath}
\usepackage{amsfonts}
\usepackage{amssymb}
\usepackage{latexsym}
\usepackage{epsfig}
\usepackage{graphicx}
\usepackage{oldgerm}
\usepackage{theorem}

\def\R{\mathbb{R}}
\def\C{\mathbb{C}}
\def\N{\mathbb{N}}
\def\Z{\mathbb{Z}}
\def\W{\mathbb{W}}
\def\P{\mathbb{P}}

\def\bE{{\bf E}}

\def\1{{\bf 1}}
\def\0{{\bf 0}}

\def\cN{{\cal N}}
\def\cH{{\cal H}}
\def\cS{{\cal S}}

\def\x{\mib{x}}
\def\y{\mib{y}}
\def\B{\mib{B}}
\def\k{\mib{k}}
\def\T{\mib{T}}
\def\a{\mib{a}}
\def\b{\mib{b}}
\def\X{\mib{X}}
\def\vnu{\mib{\nu}}
\def\vmu{\mib{\mu}}

\def\mM{\mathfrak{M}}
\def\supp{{\rm supp}\ }

\def\bK{{\bf K}}
\def\mbK{\mathbb{K}}
\def\mX{\mathfrak{X}}
\def\mY{\mathfrak{Y}}

\def\law{\stackrel{\rm (law)}{=}}

\theorembodyfont{\itshape}

\newtheorem{thm}{Theorem}[section]

\newtheorem{cor}[thm]{Corollary}
\newtheorem{prop}[thm]{Proposition}

\newcommand{\mib}[1]{\mbox{\boldmath $#1$}}
\newcommand{\SSC}[1]{\section{#1}\setcounter{equation}{0}}
\newcommand{\qed}{\hbox{\rule[-2pt]{3pt}{6pt}}}

\begin{document}

\title{\bf Reciprocal Time Relation of \\
Noncolliding Brownian Motion with Drift}
\author{
Makoto Katori
\footnote{
Department of Physics,
Faculty of Science and Engineering,
Chuo University, 
Kasuga, Bunkyo-ku, Tokyo 112-8551, Japan;
e-mail: katori@phys.chuo-u.ac.jp
}}
%%%%%%%%%%%%%%%%%%%%%%%%%%%%%%%%%%%%
\date{15 June 2012}
%%%%%%%%%%%%%%%%%%%%%%%%%%%%%%%%%%%%
\pagestyle{plain}
\maketitle
\begin{abstract}
We consider an $N$-particle system of 
noncolliding Brownian motion starting from
$x_1 \leq x_2 \leq \dots \leq x_N$
with drift coefficients $\nu_j, 1 \leq j \leq N$ 
satisfying 
$\nu_1 \leq \nu_2 \leq \dots \leq \nu_N$.
When all of the initial points are degenerated to be zero,
$x_j=0, 1 \leq j \leq N$,
the equivalence is proved between 
a dilatation with factor $1/t$ of 
this drifted process 
and the noncolliding Brownian motion starting from
$\nu_1 \leq \nu_2 \leq \dots \leq \nu_N$
without drift observed at reciprocal time $1/t$, 
for arbitrary $t > 0$.
Using this reciprocal time relation,
we study the determinantal property
of the noncolliding Brownian motion with drift
having finite and infinite numbers of particles.
\end{abstract}

%\clearpage

%%%%%%%%%%%%%%%%%%%%%%%%%%%%%%%%%%%%%%%%%%%%%%%%%%%%%%%%%%
%%%  SEC1   %%%%%%%%%%%%%%%%%%%%%%%%%%%%%%%%%%%%%%%%%%%%%%
%%%%%%%%%%%%%%%%%%%%%%%%%%%%%%%%%%%%%%%%%%%%%%%%%%%%%%%%%%
\SSC{Introduction}%%%
%%%%%%%%%%%%%%%%%%%%%%%%%%%%%%%%%%%%%%%%%%%%%%%%%%%%%%%%%%

For $N=2,3, \dots$, we consider an unbounded domain of $\R^N$,
\begin{equation}
\W_N=\{\x =(x_1, x_2, \dots, x_N) \in \R^N :
x_1 < x_2 < \cdots < x_N \},
\label{eqn:Weyl}
\end{equation}
which is called the Weyl chamber of type A$_{N-1}$
in the representation theory \cite{FH91}.
Let $\B^{\x}(t), t \geq 0$ be an $N$-dimensional
standard Brownian motion starting from 
$\x \in \W_N$. We assume that
the boundary
$\partial \W_N=\{\x \in \overline{\W}_N :
1 \leq ^{\exists}\!j \leq N-1, \quad \mbox{s.t.} \quad
x_j=x_{j+1} \}$
is made of absorbing walls in the sense that,
if $\B^{\x}(t)$ hits any point of $\partial \W_N$,
the Brownian particle is immediately absorbed
there and the Brownian motion is stopped.
For such an {\it absorbing Brownian motion
in} $\W_N$, the {\it survival probability}
$\cN_N(T, \x)$
is defined as the probability at time $T \geq 0$
such that $\B^{\x}(T), \x \in \W_N$
is not yet absorbed at $\partial \W_N$
and still moving in $\W_N$.
We can show that \cite{KT02} for any $\x \in \W_N$,
\begin{equation}
\cN_N(T, \x) \simeq c_N T^{-N(N-1)/4} h_N(\x)
\to 0 \quad \mbox{as} \quad T \to \infty,
\label{eqn:cN1}
\end{equation}
where $c_N=\pi^{-N/2} \prod_{j=1}^{N} \Gamma(j/2)/\Gamma(j)$
with the Gamma function $\Gamma$
and 
\begin{equation}
h_N(\x)=
\det_{1 \leq j, k \leq N} [x_j^{k-1}] 
=\prod_{1 \leq j < k \leq N} (x_k-x_j).
\label{eqn:Vandermonde}
\end{equation}

On the other hand, provided that
\begin{equation}
\vnu=(\nu_1, \nu_2, \dots, \nu_N) \in \W_N,
\label{eqn:nu1}
\end{equation}
the absorbing Brownian motion in $\W_N$
with drift $\vnu$, 
\begin{equation}
\widehat{\B}^{\x}(t) = \B^{\x}(t)+ \vnu t, \quad t \in [0, \infty),
\label{eqn:Bxnu}
\end{equation}
starting from $\x \in \W_N$ has a positive
probability to survive forever.
If we regard the $j$-th component of
(\ref{eqn:Bxnu}) as position of 
a Brownian motion in $\R, 1 \leq j \leq N$,
(\ref{eqn:Bxnu}) gives a configuration at time $t$
of an $N$-particle system of one-dimensional 
standard Brownian motions.
While $\widehat{\B}^{\x}(t) \in \W_N$,
any collision of particles does not occur
among $N$ particles,
and only when $\widehat{\B}^{\x}_{\vnu}(t) \in \partial \W_N$,
particle collision occurs.
If $\x \in \W_N$ and (\ref{eqn:nu1}) is satisfied, that is, 
the ordering of initial positions 
coincides with that of values of drift coefficients; 
$x_j < x_{j+1}, \nu_j < \nu_{j+1}, 1 \leq j \leq N-1$,
the probability to avoid any collision 
of particles $\cN_N^{\vnu}(T, \x)$
is positive for all $T \geq 0$.
As a matter of fact, Biane, Bougerol
and O'Connell \cite{BBO05} gave the following determinantal 
expression for the long-term limit,
\begin{equation}
\lim_{T \to \infty} 
\cN_N^{\vnu}(T, \x)=e^{-\vnu \cdot \x}
\det_{1 \leq j, k \leq N} 
[e^{\nu_j x_k} ].
\label{eqn:survival1}
\end{equation}

Conditionally on {\it staying in $\W_N$ forever},
we call the absorbing Brownian motion in $\W_N$ 
the {\it noncolliding Brownian motion} \cite{Bia95,Gra99,KT07}.
Without drift, $\vnu=0$, its transition probability
density from $\x \in \W_N$ to $\y \in \W_N$
in time duration $0 \leq t < \infty$
is given by \cite{KT02}
\begin{eqnarray}
p_N(t, \y|\x)
&=& \lim_{T \to \infty} 
\frac{\cN_N(T-t, \y)}{\cN_N(T, \x)} q_N(t, \y|\x)
\nonumber\\
&=& \frac{h_N(\y)}{h_N(\x)} q_N(t,\y|\x),
\label{eqn:pN1}
\end{eqnarray}
where (\ref{eqn:cN1}) has been used 
and $q_N$ is the transition probability density
of the absorbing Brownian motion in $\W_N$
expressed by the Karlin-McGregor determinant \cite{KM59}
\begin{equation}
q_N(t,\y|\x)=
\det_{1 \leq j, k \leq N} [p(t, y_j|x_k)].
\label{eqn:qN1}
\end{equation}
Here $p$ is the transition probability density
of the one-dimensional standard Brownian motion
({\it i.e.} the probability density of
the Gaussian distribution with variance $t$),
\begin{equation}
p(t, y|x)=
\frac{1}{\sqrt{2 \pi t}}
e^{-(y-x)^2/2t}.
\label{eqn:Gauss1}
\end{equation}
Equation (\ref{eqn:pN1}) shows the fact that
the noncolliding Brownian motion is 
Doob's $h$-transform by 
the Vandermonde determinant
(the product of differences) (\ref{eqn:Vandermonde}) 
of the absorbing Brownian motion in $\W_N$
\cite{Bia95,Gra99,KT_Sugaku_11}.

%%%%%%%%%%%%%%% Remark 1 %%%%%%%%%%%%%%%%%%%%%%%%%%%%%%%%%
\vskip 0.5cm
\noindent{\bf Remark 1} \,
An important fact is that 
(\ref{eqn:pN1}) with (\ref{eqn:qN1}) 
satisfies the partial differential equation
\begin{equation}
\left[ \frac{\partial}{\partial t}
-\left( \frac{1}{2} \Delta
+\nabla \log h_N(\x) \cdot \nabla \right) \right]
p_N(t, \y|\x)=0
\label{eqn:Kol1}
\end{equation}
with the initial condition 
$p_N(0, \y|\x)= \delta(\x-\y) \equiv \prod_{j=1}^N \delta(x_j-y_j)$,
where $\Delta=\sum_{j=1}^N \partial^2/\partial x_j^2$
and $\nabla=(\partial/\partial x_1, \dots, \partial/\partial x_N)$,
which is identified with the backward Kolmogorov equation
of Dyson's Brownian motion model with parameter $\beta=2$
\cite{Dys62,Gra99,KT04}.
That is, the noncolliding Brownian motion is equivalent to
the eigenvalue process of
Hermitian-matrix-valued Brownian motion, which has been
extensively studied as a typical example of 
log-gas system in the random matrix theory
\cite{Meh04,For10,KT_Sugaku_11}.
\vskip 0.5cm
%%%%%%%%%%%%%%%%%%%%%%%%%%%%%%%%%%%%%%%%%%%%%%%%%%%%%%%%%%%

When $\vnu \not= 0$ the transition probability density
of the noncolliding Brownian motion will be given by
\begin{equation}
p_N^{\vnu}(t, \y|\x)
=\lim_{T \to \infty} 
\frac{\cN_N^{\vnu}(T-t, \y)}{\cN_N^{\vnu}(T, \x)}
q_N^{\vnu}(t, \y|\x),
\quad \x, \y \in \W_N, \quad t \in [0, \infty)
\label{eqn:pNnu1}
\end{equation}
with the drift transform of (\ref{eqn:qN1}),
\begin{equation}
q_N^{\vnu}(t, \y|\x)
=\exp \left\{ -\frac{t}{2}|\vnu|^2
+\vnu \cdot (\y-\x) \right\}
q_N(t, \y|\x).
\label{eqn:qNnu1}
\end{equation}
If (\ref{eqn:nu1}) is satisfied, by (\ref{eqn:survival1})
Biane, Bougerol and O'Connell derived the expression \cite{BBO05}
\begin{equation}
p_N^{\vnu}(t, \y|\x)
=e^{-t|\vnu|^2/2}
\frac{\displaystyle{\det_{1 \leq j, k \leq N}[e^{\nu_j y_k}]}}
{\displaystyle{\det_{1 \leq j, k \leq N}[e^{\nu_j x_k}]}}
q_N(t, \y|\x),
\quad \x, \y, \vnu \in \W_N, \quad t \in [0, \infty).
\label{eqn:pNnu2}
\end{equation}
We note that even when some of $\nu_j$'s in 
$\vnu=(\nu_1, \dots, \nu_N)$ coincide, the ratio of determinants 
$\det_{1 \leq j,k \leq N}[e^{\nu_j y_k}]/
\det_{1 \leq j, k \leq N}[e^{\nu_j x_k}]$
can be interpreted using l'H\^opital's rule.
Then the domain of $\vnu$ is extended from (\ref{eqn:nu1}) to 
\begin{equation}
\vnu \in \overline{\W}_N
=\{ \x \in \R^N : x_1 \leq x_2 \leq \cdots \leq x_N \}.
\label{eqn:nu2}
\end{equation}
In particular, if we take the limit
$\nu_j \to 0, 1 \leq j \leq N$,
(\ref{eqn:pNnu2}) is reduced to (\ref{eqn:pN1})
(see Appendix A).
Similarly, the initial configuration $\x$ of
(\ref{eqn:pNnu2}) can be also extended to
any element of $\overline{\W}_N$.

%%%%%%%%%%%%%%% Remark 2 %%%%%%%%%%%%%%%%%%%%%%%%%%%%%%%%%
\vskip 0.5cm
\noindent{\bf Remark 2} \,
Recently O'Connell \cite{OCo12a} introduced an interacting
diffusive particle system,
which is regarded as a multivariate extension of
a one-dimensional diffusion studied by
Matsumoto and Yor \cite{MY00,MY05}.
Let $\psi_{\vnu}^{(N)}(\x), \x \in \R^N, \vnu \in \C^N$
be the class-one Whittaker function 
\cite{BO11,OCo12a,Kat11,Kat12,COSZ11,OCo12b,BC11},
whose Givental integral representation is
given by
\begin{eqnarray}
\psi_{\vnu}^{(N)}(\x)
&=&\int_{\Gamma_N(\x)} 
\exp \left[ 
\sum_{j=1}^{N} \nu_{j}
\left( \sum_{k=1}^{j} T_{j, k}
-\sum_{k=1}^{j-1} T_{j-1, k} \right)
\right. \nonumber\\
&& \quad \left.
- \sum_{j=1}^{N-1} \sum_{k=1}^j
\Big\{ e^{-(T_{j,k}-T_{j+1,k})}
+e^{-(T_{j+1, k+1}-T_{j,k})} \Big\} \right] d \T,
\label{eqn:WInt2}
\end{eqnarray}
where the integral is performed 
over the space of all real lower
triangular arrays with size $N$,
$\T=(T_{j,k}, 1 \leq k \leq j \leq N)$
conditioned $T_{N,k}=x_k, 1 \leq k \leq N$.
The transition probability density 
$p_N^{\vnu, \, a}$ of the O'Connell process
is a unique solution of a ``geometric lifting 
with parameter $a >0$" \cite{BBO09,Kat12} of (\ref{eqn:Kol1}),
\begin{equation}
\left[ \frac{\partial}{\partial t}
-\left( \frac{1}{2} \Delta
+\nabla \log \psi^{(N)}_{\vnu}(\x/a) \cdot \nabla \right) \right]
P_N^{\vnu, \, a}(t, \y|\x)=0
\label{eqn:Kol2}
\end{equation}
with the initial condition 
$P_N^{\vnu, \, a}(0, \y|\x)=\delta(\x-\y)$,
where the log-potential $\log h_N(\x)$
in (\ref{eqn:Kol1}) is replaced by 
$\log \psi^{(N)}_{\vnu}(\x/a)$.
The solution is given by
\begin{equation}
P_N^{\vnu, \, a}(t, \y|\x)
=e^{-t|\vnu|^2/2 a^2}
\frac{\psi^{(N)}_{\vnu}(\y/a)}{\psi^{(N)}_{\vnu}(\x/a)}
Q_N^{a}(t, \y|\x)
\label{eqn:pNnuxi1}
\end{equation}
with
\begin{equation}
Q_N^{a}(t, \y|\x)
=\int_{\R^N} e^{-t|\k|^2/2}
\psi^{(N)}_{i a \k}(\x/a) \psi^{(N)}_{-i a \k}(\y/a)
s_N(a \k) d \k,
\label{eqn:QN1}
\end{equation}
where $s_N(\vmu)$ is the density function of
the Sklyanin measure
\begin{equation}
s_N(\vmu)
=\frac{1}{(2 \pi)^N N!}
\prod_{1 \leq j < \ell \leq N}
|\Gamma(i(\mu_{\ell}-\mu_j))|^{-2}.
\label{eqn:sN1}
\end{equation}
We can show (see Appendix B) that
(\ref{eqn:pNnu2}) is regarded as the
``tropical analogue" of the transition probability density
of the O'Connell process in the sense that
\begin{equation}
p_N^{\vnu}(t, \y|\x)
=\lim_{a \to 0}
P_N^{a \vnu, \, a}(t, \y|\x),
\quad
\x, \vnu \in \overline{\W}_N, \,
\y \in \W, \, t \in [0, \infty).
\label{eqn:pNnu3}
\end{equation}
\vskip 0.5cm
%%%%%%%%%%%%%%%%%%%%%%%%%%%%%%%%%%%%%%%%%%%%%%%%%%

In the present paper, we study the noncolliding
Brownian motion with drift satisfying (\ref{eqn:nu2})
under the special initial condition
such that all particles are put at the origin,
$\0=(0,0, \dots, 0)$.
By taking the limit $|\x|\to 0$,
(\ref{eqn:pNnu2}) gives
\begin{eqnarray}
p_N^{\vnu}(t, \y|\0)
&=& \lim_{|\x| \to 0} p_N^{\vnu}(t, \y|\x)
\nonumber\\
&=& t^{-N} \frac{h_N(\y/t)}{h_N(\vnu)}
q_N(t^{-1}, \y/t|\vnu),
\label{eqn:pNnu0}
\end{eqnarray}
where $\y/t \equiv (y_1/t, \dots, y_N/t)$
(see Appendix A).
Then comparing with (\ref{eqn:pN1})
we obtain the equality
\begin{equation}
p_N^{\vnu}(t, \y|\0) d \y
=p_N(t^{-1}, \y/t| \vnu) d(\y/t),
\quad \y \in \W_N, 
\quad \vnu \in \overline{\W}_N, 
\quad t \in [0, \infty),
\label{eqn:equality1}
\end{equation}
where $d(\y/t)=\prod_{j=1}^N dy_j/t=t^{-N} d \y$.
It implies the equivalence between
a dilatation by factor $1/t$ of
the noncolliding Brownian motion with drift 
$\vnu \in \overline{\W}_N$ and
the process without drift observed
at reciprocal time $1/t$ (the reciprocal time relation).
Essentially the same equation with (\ref{eqn:pNnu0}) was
given by Jones and O'Connell (Proposition 2.3 in \cite{JO06})
and the equivalence (\ref{eqn:equality1})
in more generality was discussed at a special time $t=1$.
In the present paper, we extend their result to the
level of processes, in which the time change
$t \to 1/t$ is associated.
We note that, for the case $N=1$,
the statement is just a rewriting of a well-known
drift transformation (see, for instance, \cite{KS91}).
Here we discuss the reciprocal time relation
in interacting particle systems.

In an earlier paper \cite{KT10},  
the noncolliding Brownian motion with $\vnu=0$
(the Dyson model) with finite and infinite numbers of particles
were systematically studied.
If the number of particles is finite $N < \infty$
and an initial configuration is
deterministic $\xi_N=\sum_{j=1}^{N} \delta_{x_j}$,
then it is proved that the process
is {\it determinantal} and the {\it multitime
correlation kernel} is
explicitly given as a functional of $\xi_N$
(Proposition 2.1 of \cite{KT10}),
which is denoted by $\mbK^{\xi_N}$.
By the equality (\ref{eqn:equality1}),
we can conclude that at arbitrary time
$t \in [0, \infty)$, the particle distribution
of the noncolliding Brownian motion
with drift $\vnu \in \W_N$ starting from $\0$
is a {\it determinantal point process}
\cite{Sos00,ST03} with the spatial
correlation kernel 
$\mbK^{\nu_N}(t^{-1},x/t;t^{-1},y/t) \times (1/t)$
with 
$\nu_N \equiv \sum_{j=1}^N \delta_{\nu_j}$.
This spatial correlation kernel is written as
\begin{equation}
\bK_{\nu_N}(t, x; t, y)
=t \sum_{j=1}^N p(t, y|t \nu_j)
\int_{\R}d \mu' p(t, -i x|t \mu')
\prod_{1 \leq k \leq N, k \not=j}
\left(1-\frac{i \mu'-\nu_j}{\nu_k-\nu_j} \right), 
\label{eqn:bK0}
\end{equation}
$(x, y) \in \R^2, t \in [0, \infty)$, 
where $i=\sqrt{-1}$ 
and $p$ is given by (\ref{eqn:Gauss1}).

In this paper, we prove the reciprocal time relation
of the noncolliding Brownian motion with drift 
(Theorem \ref{thm:reciprocal}).
Then by a combination of this theorem
and Proposition 2.1 in \cite{KT10}, we conclude that
if the number of particles is finite
the noncolliding Brownian motion
starting from $\0$ is
determinantal for any $\vnu \in \overline{\W}_N$
and the explicit form of multitime correlation kernel
$\bK_{\nu_N}(s, \cdot; t, \cdot),
(s,t) \in [0, \infty)^2$ is determined
(Proposition \ref{thm:Finite} (i)).
There if $\vnu \in \W_N$, the expression $\bK_{\nu_N}$
is simplified (Proposition \ref{thm:Finite} (ii)).
Then by applying the theory of the Dyson model
with an infinite number of particles \cite{KT09,KT10,KT11,KT_cBM},
infinite particle limits are discussed
(Proposition \ref{thm:Infinite}, Corollary \ref{thm:LongTerm},
and Corollary \ref{thm:ShortTerm}).

We would like to put emphasis on the fact that
the reciprocal time relation will play an important role
in understanding some of recent results
for the exactly solvable interacting particle systems.
See a comment given below Corollary 4.2 in \cite{OCo12a}
with \cite{Joh04} and \cite{Kat12b}.

The paper is organized as follows.
In Sect.2 preliminaries and main results are given.
Section 3 is devoted to proofs of results.
Appendices A and B are prepared 
for explaining the facts used in Sect.1.

%%%%%%%%%%%%%%%%%%%%%%%%%%%%%%%%%%%%%%%%%%%%%%%%%%%%%%%%%%
%%%  SEC2   %%%%%%%%%%%%%%%%%%%%%%%%%%%%%%%%%%%%%%%%%%%%%%
%%%%%%%%%%%%%%%%%%%%%%%%%%%%%%%%%%%%%%%%%%%%%%%%%%%%%%%%%%
\SSC{Preliminaries and Main Results}%%%
%%%%%%%%%%%%%%%%%%%%%%%%%%%%%%%%%%%%%%%%%%%%%%%%%%%%%%%%%%

For a finite number of particles $N < \infty$,
each fixed configuration $\x=(x_1, x_2, \dots, x_N) \in \R^N$
is identified with a finite summation of delta measures,
$\xi_N(\cdot)=\sum_{j=1}^N \delta_{x_j}(\cdot)$,
if particles are indistinguishable.
For example, the configuration that all $N$ particles
are put at the origin, 
which is written as $\0$ in (\ref{eqn:pNnu0}), 
is identified with $N \delta_0(\cdot)$.
In order to discuss configurations with $N=\infty$ also,
we consider the space of nonnegative integer-valued
Radon measures on $\R$ denoted by $\mM$.
Any element $\xi$ of $\mM$ can be represented as
$\xi(\cdot)=\sum_{j \in \Lambda} \delta_{x_j}(\cdot)$
with a sequence of points in $\R$, 
$\x=(x_j)_{j \in \Lambda}$, where the index set $\Lambda$
is a finite set, the set of natural numbers $\N=\{1,2, \dots\}$ or 
that of integers $\Z=\{\dots, -1,0,1,2, \dots\}$, 
and $\xi(K)=\sharp\{j \in \Lambda: x_j \in K\} < \infty$
for any compact subset $K \subset \R$.
We call an element $\xi$ of $\mM$ a unlabeled configuration,
and a sequence $\x$ a labeled configuration.
For $\xi(\cdot)=\sum_{j\in \Lambda}
\delta_{x_j}(\cdot) \in\mM$
and $A \subset \R$, we write the
restriction of $\xi$ in $A$ as
$(\xi\cap A) (\cdot)=\sum_{j \in \Lambda : x_j \in  A}
\delta_{x_j}(\cdot)$. 
The dilatation of $\xi(\cdot)$
with factor $c >0$ is written as
$c \circ \xi(\cdot)=\sum_{j \in \Lambda} \delta_{c x_j}(\cdot)$.

First we assume that the initial configuration $\xi_N$
has a finite number of particles, $\xi_N(\R)=N \in \N$.
We consider an element of $\mM$ corresponding to
the drift coefficients $\vnu=(\nu_1, \dots, \nu_N) 
\in \overline{\W}_N$
and write it as $\nu_N(\cdot)=\sum_{j=1}^N \delta_{\nu_j}(\cdot)$.
For $M \in \N$, the multitime joint probability density
for arbitrary $M$ sequence of times
$0 < t_1 < t_2 < \cdots < t_M < \infty$ of the
noncolliding Brownian motion with drift 
$\vnu \in \overline{\W}_N$ is
given by
\begin{eqnarray}
&& p^{\xi_N}_{\nu_N}(t_1, \xi^{(1)}; \cdots : t_M, \xi^{(M)})
\nonumber\\
&& = \prod_{m=1}^{M-1} p^{\vnu}_N(t_{m+1}-t_m; \x^{(m+1)}|\x^{(m)})
p^{\vnu}_N(t_1, \x^{(1)}|\x)
\nonumber\\
&& = e^{-t_M|\vnu|^2/2}
\det_{1 \leq j, k \leq N} \left[e^{\nu_j x^{(M)}_k} \right]
\prod_{m=1}^{M-1} q_N(t_{m+1}-t_m, \x^{(m+1)}|\x^{(m)})
\frac{q_N(t_1, \x^{(1)}|\x)}
{\displaystyle{\det_{1 \leq j, k \leq N}[e^{\nu_j x_k}]}}
\label{eqn:multi1}
\end{eqnarray}
with $\xi^{(m)}=\sum_{j=1}^N \delta_{x^{(m)}_j},
\x^{(m)}=(x^{(m)}_1, \dots, x^{(m)}_N) \in \W_N,
1 \leq m \leq M, 
\x=(x_1, \dots, x_N) \in \overline{\W}_N$,
where (\ref{eqn:pNnu2}) has been used.
We can take the limit $\vnu \to \0$ of (\ref{eqn:multi1})
(using (\ref{eqn:asym1}) in Appendix A) and obtain the following expression
for the multitime joint probability density
of the noncolliding Brownian motion without drift \cite{KT07},
\begin{eqnarray}
&& p^{\xi_N}(t_1, \xi^{(1)}; \cdots; t_M, \xi^{(M)})
\nonumber\\
&& = h_N(\x^{(M)})
\prod_{m=1}^{M-1} q_N(t_{m+1}-t_m; \x^{(m+1)}|\x^{(m)})
\frac{q_N(t_1, \x^{(1)}|\x)}{h_N(\x)}.
\label{eqn:multi2}
\end{eqnarray}

The stochastic differential equation of the
noncolliding Brownian motion with drift $\vnu$,
$\X(t)=(X_1(t), \dots, X_N(t)), t \in [0, \infty)$ 
is given by
\begin{equation}
d X_j(t)=dB_j(t)+\nu_j t
+\sum_{1 \leq k \leq N, k \not=j}
\frac{dt}{X_j(t)-X_k(t)},
\quad 1 \leq j \leq N \in \N, \quad
t \in [0, \infty),
\label{eqn:SDE1}
\end{equation}
where $\{B_j(t)\}_{j=1}^N$ are independent one-dimensional
standard Brownian motions starting from 0.
Let $\Xi(t, \cdot)=\sum_{j=1}^{N} \delta_{X_j(t)}(\cdot),
t \in [0, \infty)$ and
the probability measure of the process starting from
$\xi_N \in \mM$ 
with drift $\nu_N$ be
denoted by $\P^{\xi_N}_{\nu_N}$.
The finite dimensional distribution of $\P^{\xi_N}_{\nu_N}$
is given by (\ref{eqn:multi1}).
For the noncolliding Brownian motion without drift
starting from $\xi_N$, the probability measure
is simply denoted by $\P^{\xi_N}$,
whose finite dimensional distribution is
given by (\ref{eqn:multi2}).
In general, two processes having the same state space
are said to be {\it equivalent}, if they have the same
finite-dimensional distributions, that is,
if for any finite sequence of times 
$0 < t_1 < \cdots < t_M < \infty, M \in \N$,
the multitime joint probability density functions
coincide with each other \cite{RY98}.
Here we use a symbol `$\law$' to express
equivalence of processes.

The main theorem of the present paper is the following.
%%%%%%%%%%%%%%%%%%%%%%%%%%%%%%%%%%%%%%%%%%%%%%%%%%%%%
%%%%%%%%%%%%% Theorem Reciprocal %%%%%%%%%%%%%%%%%%%%
%%%%%%%%%%%%%%%%%%%%%%%%%%%%%%%%%%%%%%%%%%%%%%%%%%%%%
\begin{thm}
\label{thm:reciprocal}
If $\nu_N \in \mM$ with $\nu_N(\R) =N \in \N$,
the following equivalence is established,
\begin{equation}
\left( \frac{1}{t} \circ \Xi(t), t \in [0, \infty), \P^{N \delta_0}_{\nu_N} \right)
\law 
\left( \Xi\left(\frac{1}{t} \right),
t \in [0, \infty), \P^{\nu_N} \right).
\label{eqn:reciprocal}
\end{equation}
\end{thm}
\vskip 0.3cm
%%%%%%%%%%%%%%%%%%%%%%%%%%%%%%%%%%%%%%%%%%%%%%%%%%%%
We call this equivalence between
a dilatation with factor $1/t$ of
the drifted process
and the process without drift
observed at reciprocal time $1/t$
the {\it reciprocal time relation}
in this paper.

For $\x^{(m)}=(x^{(m)}_1, \dots, x^{(m)}_N) \in \W_N$
and $N' \in \{ 1,2, \dots, N\}$, we put
$\x^{(m)}_{N'}=(x^{(m)}_1, \dots, x^{(m)}_{N'})
\in \W_{N'}$,
$1 \leq m \leq M$.
For a sequence $(N_m)_{m=1}^{M}$ of positive integers 
less than or equal to $N$,
we define the 
$(N_1, \dots, N_{M})$-{\it multitime correlation function} by
\begin{eqnarray}
&& \rho_{\nu_N}^{N \delta_0} (t_{1}, \x^{(1)}_{N_1} ; 
\dots; t_M, \x^{(M)}_{N_M}) 
\nonumber\\
&&=
\int_{\prod_{m=1}^{M} \R^{N-N_{m}}}
\prod_{m=1}^{M} \prod_{j=N_{m}+1}^{N} dx_{j}^{(m)}
p_{\nu_N}^{N \delta_0}(t_1, \xi_N^{(1)}; 
\dots; t_M, \xi_N^{(M)})
\prod_{m=1}^{M}
\frac{1}{(N-N_{m})!},
\label{eqn:corr}
\end{eqnarray}
which is symmetric in the sense that
$
\rho_{\nu_N}^{N \delta_0}(\dots; t_m, \sigma(\x^{(m)}_{N_m}); \dots)
=\rho_{\nu_N}^{N \delta_0}(\dots; t_m, \x^{(m)}_{N_m}; \dots)
$
with
$\sigma(\x^{(m)}_{N_m})
\equiv (x^{(m)}_{\sigma(1)}, \dots, x^{(m)}_{\sigma(N_m)})$
for any permutation $\sigma \in \cS_{N_m},
1 \leq \forall m \leq M$. 

If there is a function $\bK(s,x;t,y)$,
which is continuous with respect to
$(x,y) \in \R^2$ for any fixed $(s,t) \in [0, \infty)^2$,
such that it determines 
the finite dimensional distributions of the process 
by giving determinantal expressions to 
multitime correlation functions $\rho$ as
$$
\rho (t_1,\x^{(1)}_{N_1}; \dots;t_M,\x^{(M)}_{N_M} ) 
=\det_{
\substack{1 \leq j \leq N_{m}, 1 \leq k \leq N_{n} \\ 1 \leq m, n \leq M}
}
\Big[
\bK(t_m, x_{j}^{(m)}; t_n, x_{k}^{(n)} )
\Big]
$$
for any integer $M \geq 1$, 
any sequence $(N_m)_{m=1}^{M}$ of positive integers, and
any time sequence $0 < t_1 < \cdots < t_M < \infty$,
then the process is said to be 
{\it determinantal with the correlation kernel} $\bK$ \cite{KT10}.

We use the convention such that
\begin{equation}
\prod_{x\in\xi}f(x) =\exp
\left\{\int_\R \xi(dx) \log f(x) \right\}
=\prod_{x \in \supp \xi}f(x)^{\xi(\{x\})}
\label{eqn:convention}
\end{equation}
for $\xi\in \mM$ and a function $f$ on $\R$,
where $\supp \xi = \{x \in \R : \xi(\{x\}) > 0\}$.
Let $\1(\omega)$ be the indicator function of 
a condition $\omega$;
$\1(\omega)=1$ if $\omega$ is satisfied and
$\1(\omega)=0$ otherwise.

We consider a set of configurations
with no multiple points, 
\begin{equation}
\mM_0= \Big\{ \xi\in\mM : \xi(\{x\})\le 1 \mbox { for any }  x\in\R
\Big\}.
\label{eqn:mM0}
\end{equation}
Since any element $\xi$ of $\mM_0$ is determined uniquely 
by its support, 
it is identified with a countable subset $\{x_j\}_{j\in\Lambda}$ of $\R$.
For $\xi_N \in \mM_0, a \in \C$,
we introduce an entire function of $z \in \C$
\begin{equation}
\Phi(\xi_N, a, z)= 
\prod_{u \in \xi_N \cap \{a\}^{\rm c}}
\left( 1 - \frac{z-a}{u-a} \right),
\label{eqn:entire1}
\end{equation}
where $A^{\rm c}$ denotes a complementary set of $A$.
The zero set of $\Phi(\xi_N, a, z)$ is 
$\supp (\xi_N \cap \{a\}^{\rm c})$ 
(see, for instance, \cite{L96}).
Then, as an application of
Proposition 2.1 in \cite{KT10}, Theorem \ref{thm:reciprocal}
gives the following.

%%%%%%%%%%%%%%%%%%%%%%%%%%%%%%%%%%%%%%%%%%%%%%%%%%%%%%%%%%%%
%%%%%%%%%%%%%%%%%%% Theorem Finite %%%%%%%%%%%%%%%%%%%%%%
%%%%%%%%%%%%%%%%%%%%%%%%%%%%%%%%%%%%%%%%%%%%%%%%%%%%%%%%%%%%%
\begin{prop}
\label{thm:Finite}
\noindent (i) 
The noncolliding Brownian motion 
with drift $\nu_N\in \mM$ with
$\nu_N(\R) = N \in \N$,
$(\Xi(t), t \in [0, \infty), \P^{N \delta_0}_{\nu_N})$, 
is determinantal with the correlation kernel 
\begin{eqnarray}
\bK_{\nu_N}(s, x; t, y)
&=& \frac{\sqrt{s t}}{2 \pi i} 
\oint_{\Gamma(\nu_N)} d \mu \, p(s, x|s \mu)
\int_{\R} d \mu' \, p(t, -iy|t \mu') 
%\nonumber\\
%&& \hskip 2cm \times
\frac{1}{i \mu'- \mu}
\prod_{u \in \nu_N}
\left( 1- \frac{i \mu'-\mu}{u- \mu} \right)
\nonumber\\
&& - \1(s < t) 
p \left(t-s, \left. \sqrt{\frac{s}{t}} y \right|
\sqrt{\frac{t}{s}} x \right),
\quad
(s,t)\in [0, \infty)^2, (x,y) \in \R^2,
\label{eqn:KN1a}
\end{eqnarray}
where $\Gamma(\nu_N)$ is a closed contour on the
complex plane $\C$ encircling the points in 
$\supp \nu_N$ on the real line $\R$
once in the positive direction. \\
(ii) If $\nu_N\in \mM_0$ with $\nu_N(\R) = N \in \N$,
the correlation kernel (\ref{eqn:KN1a}) is simplified as
\begin{eqnarray}
&& \bK_{\nu_N}(s, x; t, y)
= \sqrt{st} \sum_{j=1}^N p(s, x|s \nu_j)
\int_{\R} d \mu' \, p(t, -iy|t \mu')
\Phi(\nu_N, \nu_j, i \mu')
\nonumber\\
&&\qquad - \1(s < t) 
p \left(t-s, \left. \sqrt{\frac{s}{t}} y \right|
\sqrt{\frac{t}{s}} x \right), 
\quad
(s,t)\in [0, \infty)^2, (x,y) \in \R^2.
\label{eqn:KN1}
\end{eqnarray}
\end{prop}
\vskip 0.3cm
%%%%%%%%%%%%%%%%%%%%%%%%%%%%%%%%%%%%%%%%%%%%%%%%%%%%%%%%%%%%%
%%%%%%%%%%%%%%%%%%%%%%%%%%%%%%%%%%%%%%%%%%%%%%%%%%%%%%%%%%%%%

Note that, if $\nu_N \in \mM_0$, 
\begin{equation}
\Phi(\nu_N, \nu_j, i \mu')
=\prod_{1 \leq k \leq N, k \not=j}
\left( 1- \frac{i \mu'-\nu_j}{\nu_k-\nu_j} \right),
\label{eqn:entire2}
\end{equation}
and thus, when we set $s=t$, (\ref{eqn:KN1})
is reduced to be (\ref{eqn:bK0}).

%%%%%%%%%%%%%%% Remark 3 %%%%%%%%%%%%%%%%%%%%%%%%%%%%%%%%%
\vskip 0.5cm
\noindent{\bf Remark 3} \,
The reciprocal time relation (\ref{eqn:reciprocal})
will define a kind of duality between two processes
$(\Xi(t), t \in [0, \infty), \P^{\xi}_{\nu})$
and $(\Xi(1/t), t \in [0, \infty), \P^{\nu}_{\xi})$,
where the initial configuration and drift coefficients
are exchanged.
A self-dual process, which satisfies the equality,
\begin{equation}
\left( \frac{1}{t} \circ \Xi(t), t \in [0, \infty),
\P^{N \delta_0} \right)
\law 
\left( \Xi \left(\frac{1}{t} \right), t \in [0, \infty),
\P^{N \delta_0} \right), 
\label{eqn:self_dual}
\end{equation}
$N \in \N$, 
is known as the determinantal process
with the extended Hermite-function kernel
(see, for example, Eq.(4.25) of \cite{KT07},
Section 5.1 of \cite{KT_Sugaku_11}).
\vskip 0.5cm
%%%%%%%%%%%%%%%%%%%%%%%%%%%%%%%%%%%%%%%%%%%%%%%%%%%%%%%%%

For $\nu \in \mM$, when 
$\bK_{\nu \cap [-L, L]}$
converges to a continuous function as $L \to \infty$,
the limit is written as $\bK_{\nu}$. If 
the probability measure $\P_{\nu \cap [-L, L]}$ of the process
determined by $\bK_{\nu \cap [-L, L]}$ converges
to a probability measure $\P_{\nu}$
on $\mM^{[0, \infty)}$,
which is determinantal with the correlation kernel
$\bK_{\nu}$, weakly in the sense of finite dimensional
distributions as $L \to \infty$
in the vague topology, we say that the obtained process 
is {\it well defined with the correlation kernel}
$\bK_{\nu}$ \cite{KT09,KT10,KT11}.
The regularity of the sample paths 
can be discussed as \cite{KT_cBM}.
In the case $\nu(\R)=\infty$, the process
obtained by this limit 
has an infinite number of particles with drift $\nu$ .

In \cite{KT10}, two sets of configurations
with finite and infinite numbers of particles were
introduced and denoted by $\mX$ and $\mY$.
By the argument given there,
the following statements can be proved
for the noncolliding Brownian motion
with drift $\nu$, which can be constructed from
an infinite number of particles
all starting from the origin.
This infinite particle process
will be denoted by $(\Xi(t), t \in [0, \infty), \P_{\nu}^0)$.

%%%%%%%%%%%%%%%%%%%%%%%%%%%%%%%%%%%%%%%%%%%%%%%%%%%%%%%%%%
%%%%%%%%%%%%%%%%%%% Theorem Infinite %%%%%%%%%%%%%%%%%%%%%%
%%%%%%%%%%%%%%%%%%%%%%%%%%%%%%%%%%%%%%%%%%%%%%%%%%%%%%%%%%%%%
\begin{prop}
\label{thm:Infinite}
\noindent (i) \,
The infinite particle system 
of noncolliding Brownian motion with
all particles starting from the origin,
$(\Xi(t), t \in [0, \infty), \P_{\nu}^0)$,
is well defined, if the drift coefficients
$\nu \in \mX \cap \mM_0$. 
The correlation kernel is given by
\begin{eqnarray}
&& \bK_{\nu}(s, x; t, y)
= \sqrt{st} \int_{\R} \nu(d \mu) \, p(s, x|s \mu)
\int_{\R} d \mu' \, p(t, -iy|t \mu')
\Phi(\nu_N, \mu, i \mu')
\nonumber\\
&& \qquad - \1(s < t) 
p \left(s-t, \left. \sqrt{\frac{s}{t}} y \right|
\sqrt{\frac{t}{s}} x \right),
\quad (s,t) \in [0, \infty)^2,
(x,y) \in \R^2.
\label{eqn:KN2}
\end{eqnarray}
\noindent (ii) \,
If $\nu \in \mY$, 
$(\Xi(t), t \in [0, \infty), \P_{\nu}^0)$ is well defined.
In particular, when $\nu \in \mY \cap \mM_0$,
the correlation kernel is given by (\ref{eqn:KN2}).
\end{prop}
\vskip 0.3cm
%%%%%%%%%%%%%%%%%%%%%%%%%%%%%%%%%%%%%%%%%%%%%%%%%%

As consequences of Theorem \ref{thm:reciprocal}
and Proposition \ref{thm:Infinite}, we can conclude
the following facts for the noncolliding Brownian motion
with drift, whose number of particles is infinite.
The reciprocal time relation (\ref{eqn:reciprocal})
implies that, the long-term limit $t \to \infty$
of $((1/t) \circ \Xi(t), \P^{N \delta_0}_{\nu_N})$
is given by a $1/t \to 0$ limit of
$(\Xi(1/t), \P^{\nu_N})$, 
for any $N \in \N$.
Since $\lim_{t \to \infty} \Xi(1/t)=\nu_N$
in $\P^{\nu_N}$, we have the following.

%%%%%%%%%%%%%%%%%%%%%%%%%%%%%%%%%%%%%%%%%%%%%%%%%%%%%%%%%%
%%%%%%%%%%%%%%%%%%% Corollary Long-Term %%%%%%%%%%%%%%%%%%%%%%
%%%%%%%%%%%%%%%%%%%%%%%%%%%%%%%%%%%%%%%%%%%%%%%%%%%%%%%%%%%%%
\begin{cor}
\label{thm:LongTerm}
If $\nu \in \mX$ or $\nu \in \mY$,
\begin{equation}
\lim_{t \to \infty} \frac{1}{t} \circ \Xi(t)=\nu
\quad \mbox{in $\P^{0}_{\nu}$}.
\label{eqn:longterm}
\end{equation}
\end{cor}
\vskip 0.3cm
%%%%%%%%%%%%%%%%%%%%%%%%%%%%%%%%%%%%%%%%%%%%%%%%%%%%%%%%

Let $\nu_{\Z}(\cdot)=\sum_{j \in \Z} \delta_j(\cdot)$.
In \cite{KT10}, it is shown that 
$\nu_{\Z} \in \mX \cap \mM_0$ and the 
noncolliding Brownian motion without drift
(the Dyson model) starting from $\nu_{\Z}$ is well
defined.
Moreover, it is proved that this infinite particle
system shows a {\it relaxation phenomenon}
to the equilibrium and spatially homogeneous 
determinantal process
with the {\it extended sine kernel} with density 1,
\begin{eqnarray}
{\bf K}_{\sin}(t-s, y-x) &=&
\frac{1}{2 \pi} \int_{|k| \leq \pi} dk \,
e^{k^2(t-s)/2 + i k (y-x)}
- {\bf 1}(s>t) p(s-t, x|y) \nonumber\\
&=& \left\{ \begin{array}{ll} 
\displaystyle{
\int_{0}^{1} du \, e^{\pi^2 u^2 (t-s)/2} 
\cos \{ \pi u (y-x)\} }
& \mbox{if $t>s $} \cr
& \cr
\displaystyle{\frac{\sin \pi (y-x)}{\pi (y-x)}}
& \mbox{if $t=s$} \cr
& \cr
\displaystyle{
- \int_{1}^{\infty} du \, 
e^{\pi^2 u^2 (t-s)/2} \cos \{ \pi u (y-x) \} }
& \mbox{if $t<s$},
\end{array} \right.
\label{eqn:sine-kernel}
\end{eqnarray}
$(s, t) \in [0, \infty)^2, (x, y) \in \R^2$.
Let $\vartheta_3$ be a version of the Jacobi theta function
defined by
\begin{equation}
\vartheta_3(v, \tau) 
= \sum_{n \in \Z} e^{2 \pi i v n+\pi i \tau n^2}.
\quad \Im \tau > 0.
\label{eqn:theta1}
\end{equation}
Then by the reciprocal time relation 
(\ref{eqn:reciprocal}), we have the following.
%%%%%%%%%%%%%%%%%%%%%%%%%%%%%%%%%%%%%%%%%%%%%%%%%%%%%%%%%%
%%%%%%%%%%%%%%%%%%% Corollary Short-Term %%%%%%%%%%%%%%%%%%%%%%
%%%%%%%%%%%%%%%%%%%%%%%%%%%%%%%%%%%%%%%%%%%%%%%%%%%%%%%%%%%%%
\begin{cor}
\label{thm:ShortTerm}
The noncolliding Brownian motion
with an infinite number of particles, 
$(\Xi(t), t \in [0, \infty), \P^{0}_{\nu_Z})$, 
is well defined with the correlation kernel
\begin{eqnarray}
&& \bK_{\nu_{\Z}}(s,x;t,y)
=\frac{1}{2 \pi} \int_{|w| \leq \pi/\sqrt{st}} d w \,
e^{w^2(t-s)/2+iw(y \sqrt{s/t}- x \sqrt{t/s})}
\vartheta_3 \left( \frac{x}{s}-i w \sqrt{\frac{t}{s}},
\frac{2 \pi i}{s} \right),
\nonumber\\
&& \hskip 8cm
(s,t) \in [0, \infty)^2, (x,y) \in \R^2.
\label{eqn:Jacobi2}
\end{eqnarray}
This kernel has the following asymptotics, 
\begin{equation}
\lim_{u \to \infty}
u \bK_{\nu_{\Z}}\left(
\frac{1}{u+s}, \frac{x}{u+s}; 
\frac{1}{u+t}, \frac{y}{u+t} \right)
=\bK_{\sin}(t-s, y-x).
\label{eqn:short_term}
\end{equation}
\end{cor}
\vskip 0.3cm
%%%%%%%%%%%%%%%%%%%%%%%%%%%%%%%%%%%%%%%%%%%%%%

%%%%%%%%%%%%%%%%%%%%%%%%%%%%%%%%%%%%%%%%%%%%%%%%%%%%%%%%%%
%%%  SEC3   %%%%%%%%%%%%%%%%%%%%%%%%%%%%%%%%%%%%%%%%%%%%%%
%%%%%%%%%%%%%%%%%%%%%%%%%%%%%%%%%%%%%%%%%%%%%%%%%%%%%%%%%%
\SSC{Proofs of Theorems}%%%
%%%%%%%%%%%%%%%%%%%%%%%%%%%%%%%%%%%%%%%%%%%%%%%%%%%%%%%%%%

\noindent{\it Proof of Theorem \ref{thm:reciprocal}} \,
By (\ref{eqn:pNnu2}), (\ref{eqn:multi1}) is written as
\begin{eqnarray}
&& p^{\xi_N}_{\nu_N}(t_1, \xi^{(1)}; \cdots ; t_M, \xi^{(M)})
\nonumber\\
&& =e^{-t_M|\vnu|^2/2} 
\frac{\displaystyle{\det_{1 \leq j, k \leq N}
[e^{\nu_j x^{(M)}_k}]}}
{\displaystyle{\det_{1 \leq j, k \leq N}
[e^{\nu_j x^{(1)}_k}]}}
e^{t_1|\vnu|^2/2} 
\prod_{m=1}^{M-1} q_N(t_{m+1}-t_m; \x^{(m+1)}|\x^{(m)})
p^{\vnu}_N(t_1, \x^{(1)}|\x)
\nonumber\\
&& = \frac{(2\pi/t_M)^{N/2}}
{\displaystyle{\det_{1 \leq j, k \leq N}
[e^{\nu_j x^{(1)}_k}]}}
q_N(t_M^{-1}, \x^{(M)}/t_M|\vnu)
e^{|\x^{(M)}|^2/2t_M}
e^{t_1|\vnu|^2/2} \nonumber\\
&& \qquad \qquad \times 
\prod_{m=1}^{M-1} q_N(t_{m+1}-t_m, \x^{(m+1)}|\x^{(m)})
p^{\vnu}_N(t_1, \x^{(1)}|\x).
\label{eqn:Eq1}
\end{eqnarray}
We rewrite the equality (\ref{eqn:equality1}) as
\begin{eqnarray}
p_N^{\vnu}(t, \y|\0) d \y
&=& p_N(t^{-1}, \y/t|\vnu) d (\y/t)
\nonumber\\
&=& \frac{h_N(\y/t)}{h_N(\vnu)}
e^{-t|\vnu|^2/2} e^{-|\y|^2/2t} 
\frac{\displaystyle{\det_{1 \leq j, k \leq N}
[e^{\nu_j y_k}]}}{(2 \pi t)^{N/2}} d \y.
\label{eqn:Eq2}
\end{eqnarray}
Then, we can confirm that for $\xi_N=N \delta_0$, 
(\ref{eqn:Eq1}) gives
\begin{eqnarray}
&& p^{N \delta_0}_{\nu_N}(t_1, \xi^{(1)}; \cdots; t_M, \xi^{(M)})
= (t_M t_1)^{-N/2}
\frac{q_N(t_M^{-1}, \x^{(M)}/t_M |\vnu)}
{h_N(\vnu)}
\nonumber\\
&& \qquad \qquad \times e^{|\x^{(M)}|^2/2t_M}
\prod_{m=1}^{M-1} q_N(t_{m+1}-t_m, \x^{(m+1)}|\x^{(m)})
e^{-|\x^{(1)}|^2/2t_1}
h_N(\x^{(1)}/t_1).
\label{eqn:Eq3}
\end{eqnarray}
We find that the following equalities hold,
\begin{eqnarray}
&& q_N(t_{m+1}-t_m, \x^{(m+1)}|\x^{(m)})
\nonumber\\
&& \quad =
(t_{m+1} t_m)^{-N/2} 
e^{-|\x^{(m+1)}|^2/2t_{m+1}}
q_N \left(\left. t_m^{-1}-t_{m+1}^{-1}, \frac{\x^{(m)}}{t_m} 
 \right| \frac{\x^{(m+1)}}{t_{m+1}} \right)
e^{|\x^{(m)}|^2/2t_m},
\nonumber\\
&& \hskip 10cm 1 \leq m \leq M.
\label{eqn:equality3}
\end{eqnarray}
Then (\ref{eqn:Eq3}) is equal to
\begin{equation}
 h_N(\x^{(1)}/t_1)
\prod_{m=1}^{M-1} 
q_N \left(\left. t_m^{-1}-t_{m+1}^{-1}, \frac{\x^{(m)}}{t_m} 
 \right| \frac{\x^{(m+1)}}{t_{m+1}} \right)
 \frac{q_N(t_M^{-1}, \x^{(M)}/t_M|\vnu)}{h_N(\vnu)}
\prod_{m=1}^M t_m^{-N}.
\label{eqn:Eq4}
\end{equation}
Comparing this result with (\ref{eqn:multi2}), 
we obtain the equality
\begin{eqnarray}
&& p^{N \delta_0}_{\nu_N}(t_1, \xi^{(1)}; \cdots; t_M, \xi^{(M)})
\prod_{m=1}^{M} d \x^{(m)}
\nonumber\\
&& \qquad =
p^{\nu_N}
\left( t_M^{-1}, \frac{1}{t_M} \circ \xi^{(M)};
\cdots; t_1^{-1}, \frac{1}{t_1} \circ \xi^{(1)} \right)
\prod_{m=1}^{M} \left( \frac{d \x^{(m)}}{t_m} \right).
\label{eqn:Eq5}
\end{eqnarray}
Then the equivalence (\ref{eqn:reciprocal}) is proved. \qed

%%%%%%%%%%%%%%%%%%%%%%%%%%%%%%%%%%%%%%%%%%%%%%%%%%
\vskip 0.3cm
\noindent{\it Proof of Proposition \ref{thm:Finite}} \,
(i) By combination of
Theorem \ref{thm:reciprocal} and 
Proposition 2.1 in \cite{KT10}, the correlation kernel
is given by
\begin{eqnarray}
&& \bK_{\nu_N}(s, x; t, y)
= \mbK^{\nu_N}(s^{-1}, x/s; t^{-1}, y/t) \frac{1}{\sqrt{st}}
\nonumber\\
&& = \frac{1}{2 \pi i} \oint_{\Gamma(\nu_N)} d \mu \,
p(s^{-1}, x/s|\mu) \int_{\R} d \mu' \, p(t^{-1}, -i y/t|\mu')
\frac{1}{i\mu'-\mu} 
\prod_{u \in \nu_N} \left( 1- \frac{i \mu'-\mu}{u-\mu} \right)
\frac{1}{\sqrt{st}}
\nonumber\\
&& \qquad -\1(s^{-1}> t^{-1}) 
p(s^{-1}-t^{-1}, x/s|y/t) \frac{1}{\sqrt{st}}.
\label{eqn:KK1}
\end{eqnarray}
Here by the cyclic property of determinant,
a factor $1/\sqrt{st}$ should be put on the
correlation kernel for the transform of variables,
$(x,y) \mapsto (x/s,y/t)$.
By (\ref{eqn:Gauss1}), we see
\begin{eqnarray}
&& p(s^{-1}, x/s|\mu)=s p(s, x|s \mu),
\nonumber\\
&& p(t^{-1}, -i y/t|\mu')=t p(t, -i y|t\mu'),
\nonumber
\end{eqnarray}
and
$$
\1(s^{-1}>t^{-1}) p(s^{-1}-t^{-1}, x/s|y/t)
=\1 (s<t) \sqrt{st} 
p\left( t-s, \left. \sqrt{\frac{s}{t}} y \right|
\sqrt{\frac{t}{s}} x \right).
$$
Then (\ref{eqn:KN1a}) is obtained. \\
(ii) By performing the Cauchy integral
in (\ref{eqn:KN1a}) on the contour $\Gamma(\nu^N)$,
(\ref{eqn:KN1}) is obtained. \qed
\vskip 0.3cm
%%%%%%%%%%%%%%%%%%%%%%%%%%%%%%%%%%%%%%%%%%%%%%%

\noindent{\it Proof of Corollary \ref{thm:LongTerm}} \,
By Theorem \ref{thm:reciprocal},
for $u, s, t \in [0, \infty), x, y \in \R$,
\begin{equation}
\bK_{\nu}(u+s, ux; u+t, uy)
=\mbK^{\nu}\left( \frac{1}{u+s}, \frac{ux}{u+s};
\frac{1}{u+t}, \frac{uy}{u+t} \right)
\sqrt{\frac{u}{u+s} \frac{u}{u+t}}.
\label{eqn:LTp1}
\end{equation}
And 
\begin{eqnarray}
&& 
\mbK^{\nu}\left( \frac{1}{u+t}, \frac{ux}{u+t};
\frac{1}{u+t}, \frac{uy}{u+t} \right)
\mbK^{\nu}\left( \frac{1}{u+s}, \frac{uy}{u+s};
\frac{1}{u+t}, \frac{ux}{u+t} \right) 
\left( \frac{u}{u+t} \right) dx dy
\nonumber\\
&& \qquad \qquad \to \nu(dx) \1(x=y)
\quad \mbox{as $u \to \infty$ 
in the vague topology}.
\nonumber
\end{eqnarray}
Then the statement is obtained. \qed

\vskip 0.3cm
%%%%%%%%%%%%%%%%%%%%%%%%%%%%%%%%%%%%%%%%%%%%%%%

\noindent{\it Proof of Corollary \ref{thm:ShortTerm}} \,
The expression (\ref{eqn:Jacobi2}) is readily obtained
from the formula (1.5) in \cite{KT10}
following the reciprocal time relation (\ref{eqn:reciprocal}).
We can check that
\begin{eqnarray}
\lim_{u \to \infty} u
\bK_{\nu_{\Z}} \left( \frac{1}{u+s}, \frac{x}{u+s};
\frac{1}{u+t}, \frac{y}{u+t} \right)
&=& \lim_{u \to \infty}
\mbK^{\nu_{\Z}}(u+s, x; u+t, y)
\frac{u}{\sqrt{(u+s)(u+t)}}
\nonumber\\
&=& \bK_{\sin}(t-s, y-x).
\nonumber
\end{eqnarray}
Then the statement is concluded. \qed

%%%%%%%%%%%%%%%%%%%%%%%%%%%%%%%%%%%%%%%%%%%%%%%%%%%%%%%%%
%\vskip 1cm
\clearpage
%%%%%%%%%%%%%%%%%%%%%%%%%%%%%%%%%%%%%%%%%%%%%%%%%%%%%%%%%
%%%%%%%%%%%% Appendix %%%%%%%%%%%%%%%%%%%%%%%%%%%%%%%%%%%%%%%
%%%%%%%%%%%%%%%%%%%%%%%%%%%%%%%%%%%%%%%%%%%%%%%%%%%%%%%%%
\appendix
\begin{LARGE}
{\bf Appendices}
\end{LARGE}

%%%%%%%%%%%%%%%%%%%%%%%%%%%%%%%%%%%%%%%%%%%%%%%%%%%%%%%%
\SSC{Asymptotics of determinants}
%%%%%%%%%%%%%%%%%%%%%%%%%%%%%%%%%%%%%%%%%%%%%%%%%%%%%%%%

By the Schur function expansion, we can prove that
\cite{KT04,KT07} for $\b \in \C^N$
\begin{equation}
\det_{1 \leq j, k \leq N}
[e^{a_j b_k}]
=\frac{h_N(\a) h_N(\b)}{\prod_{j=1}^N \Gamma(j)}
\times \{1+{\cal O}(|\a|) \}
\quad \mbox{as} \quad |\a| \to 0.
\label{eqn:asym1}
\end{equation}
Then 
\begin{equation}
\lim_{|\vnu| \to 0}
\frac{\displaystyle{
\det_{1 \leq j, k \leq N}[e^{\nu_j y_k}]}}
{\displaystyle{
\det_{1 \leq j, k \leq N}[e^{\nu_j x_k}]}}
=\frac{h_N(\y)}{h_N(\x)},
\label{eqn:asym2}
\end{equation}
which implies the fact that
(\ref{eqn:pN1}) is reduced from (\ref{eqn:pNnu2})
by taking $\nu_j \to 0, 1 \leq j \leq N$.

Similarly, we can see that
$$
\lim_{|\x| \to 0}
\frac{\displaystyle{ \det_{1 \leq j, k \leq N}
[e^{x_j y_k/t}]}}
{\displaystyle{ \det_{1 \leq j, k \leq N}
[e^{\nu_j x_k}]}}
= \lim_{|\x| \to 0}
\frac{h_{N}(\x/\sqrt{t}) h_N(\y/\sqrt{t})}
{h_N(\vnu) h_N(\x)}
=\frac{h_N(\y/t)}{h_N(\vnu)}.
$$
Since $q_N(t, \y|\x)$ given by (\ref{eqn:qN1})
is equal to
$(2 \pi t)^{-N/2} e^{-(|\x|^2+|\y|^2)/2t}
\det_{1 \leq j, k \leq N} [e^{x_j y_k/t}]$,
(\ref{eqn:pNnu2}) gives
$$
\lim_{|\x| \to 0} p_N^{\vnu}(t, \y|\x)
=e^{-t|\vnu|^2/2} \det_{1 \leq j, k \leq N}[e^{\nu_j y_k}]
(2\pi t)^{-N/2} e^{-|\y|^2/2t} 
\frac{h_N(\y/t)}{h_N(\vnu)},
$$
which is equal to (\ref{eqn:pNnu0}).

%%%%%%%%%%%%%%%%%%%%%%%%%%%%%%%%%%%%%%%%%%%%%%%%%%%%%%%%
\SSC{On the O'Connell process with drift $\vnu$}
%%%%%%%%%%%%%%%%%%%%%%%%%%%%%%%%%%%%%%%%%%%%%%%%%%%%%%%%
Since the formula (\ref{eqn:pNnuxi1}) is not found
in \cite{OCo12a}, here we give explanation for it
and then prove (\ref{eqn:pNnu3}).
Note that the derivation of (\ref{eqn:pNnuxi1})
with $\vnu=0$ was given in \cite{Kat11,Kat12}.

For $a >0, \vnu \in \R^N$, we consider the
following partial differential equation
\begin{equation}
\left[ \frac{\partial}{\partial t}
+\cH_N-\frac{1}{a} \vnu \cdot \nabla 
\right] u_a^{\vnu}(t, \x)=0,
\quad \x \in \R^N,
\quad t \in [0, \infty),
\label{eqn:Toda1}
\end{equation}
where
\begin{equation}
\cH_N=-\frac{1}{2} \Delta
+\frac{1}{a^2} \sum_{j=1}^{N-1} 
e^{-(x_{j+1}-x_j)/a}
\label{eqn:Toda2}
\end{equation}
is identified with the Hamiltonian
of an open quantum Toda lattice \cite{OCo12a}.
Assume that $\x, \vnu \in \W_N$.
Then by the Feynman-Kac formula 
(see, for instance, \cite{KS91})
the stationary solution $u^{\vnu}_{a}(\x)$ of (\ref{eqn:Toda1})
with (\ref{eqn:Toda2}), which is uniquely determined
by imposing the condition
$\lim_{|\x| \to \infty, \x \in \W_N} u_{a}^{\vnu}(\x)=1$, 
$\vnu \in \W_N$, is given by
\begin{equation}
u_{a}^{\vnu}(\x)=\bE^{\x} \left[
\exp \left( - \frac{1}{a^2}
\sum_{j=1}^{N-1} \int_0^{\infty}
e^{-\{\widehat{B}_{j+1}(s)-\widehat{B}_j(s)\}/a} ds \right)
\right],
\label{eqn:u1}
\end{equation}
where $\bE^{\x}[ \, \cdot \,]$
denotes the expectation with respect to
the Brownian motion (\ref{eqn:Bxnu})
with drift $\vnu \in \W_N$ and starting from
$\x \in \W_N$.
We note that for $0 \leq T < \infty$
\begin{equation}
\cN_N^{\vnu, \, a}(T, \x)
=\bE^{\x}
 \left[
\exp \left( - \frac{1}{a^2}
\sum_{j=1}^{N-1} \int_0^{T}
e^{-\{\widehat{B}_{j+1}(s)-\widehat{B}_j(s)\}/a} ds \right)
\right]
\label{eqn:Nnuxi}
\end{equation}
expresses the probability that,
in the mutually killing $N$-particle system 
with the killing term
$-(1/a^2) \sum_{j=1}^{N-1} e^{-(x_{j+1}-x_j)/a}$
\cite{Kat12},
all $N$ Brownian particles with drifts
$\{\widehat{B}_j(t)\}_{j=1}^N$ starting from $\x \in \W_N$
survive up to time $T$
and that (\ref{eqn:u1}) is its long-term limit,
$\lim_{T \to \infty} \cN_N^{\vnu, \, a}(T, \x)$,
$\x, \vnu \in \W_N$.
On the other hand, the class-one Whittaker function 
$\psi^{(N)}_{\vnu}(\x)$ given by (\ref{eqn:WInt2})
is an eigenfunction of the Hamiltonian (\ref{eqn:Toda2})
with the eigenvalue $-\sum_{j=1}^N \nu_j^2/2$,
$\vnu \in \C^N$. 
By the method of separation of variables,
we can show that $u_{a}(\x)$ is also expressed
by using $\psi^{(N)}_{\vnu}(\x/a)$.
Then the equality
\begin{equation}
\lim_{T \to \infty} \cN_N^{\vnu, \, a}(T, \x)
=c_N(\vnu) e^{-\vnu \cdot \x/a} 
\psi^{(N)}_{\vnu}(\x/a),
\quad \x, \vnu \in \W_N,
\label{eqn:equality2}
\end{equation}
is established, where
$c_N(\vnu)=\prod_{1 \leq j < k \leq N}
\{\sin \pi (\nu_k-\nu_j)\}/\pi$ 
\cite{OCo12b}.

We find that $Q_N^{a}(t, \y|\x)$ given by
(\ref{eqn:QN1}) solves (\ref{eqn:Toda1})
with $\vnu=0$ and $Q_N^{a}(0,\y|\x)=\delta(\x-\y)$,
and then we can regard it as a transition probability
density of the mutually killing Brownian motions
with duration $t \in [0, \infty)$ from $\x \in \R^N$
to $\y \in \R^N$, preserving the number of particles \cite{Kat11,Kat12}.
Then the transition probability density of the
mutually killing Brownian motions with 
drift $\vnu$ {\it conditioned that all particle survive}
is given by
\begin{equation}
P_N^{\vnu,\, a}(t, \y|\x)
=\lim_{T \to \infty} \frac{\cN_N^{\vnu, \, a}(T-t, \y)}
{\cN_N^{\vnu, \, a}(T, \x)}
Q_N^{\vnu, \, a}(t, \y|\x)
\label{eqn:PNA1}
\end{equation}
with the drift transform of (\ref{eqn:QN1}) with parameter $a>0$
\begin{equation}
Q_N^{\vnu, \, a}(t, \y|\x)
= \exp \left\{ - \frac{t|\vnu|^2}{2 a^2}
+\frac{\vnu}{a} \cdot (\y-\x) \right\}
Q_N^{a}(t,\y|\x).
\label{eqn:QNA1}
\end{equation}
If $\x, \y, \vnu \in \W_N$, by (\ref{eqn:equality2}),
(\ref{eqn:PNA1}) is equal to (\ref{eqn:pNnuxi1}),
which should solve (\ref{eqn:Kol2})
with the initial condition 
$P_N^{\vnu, \, a}(0, \y|\x)=\delta(\x-\y)$.

The class-one Whittaker function has the
alternating sum formula \cite{BO11},
\begin{equation}
\psi^{(N)}_{\vnu}(\x)=c_N(\vnu)^{-1}
\sum_{\sigma \in \cS_N} {\rm sgn}(\sigma)
m^{(N)}(\x, \sigma(\vnu)),
\label{eqn:tro1}
\end{equation}
where
$\sigma(\vnu)=(\nu_{\sigma(1)}, \dots, \nu_{\sigma(N)})$
for each permutation $\sigma \in \cS_N$.
Here $m^{(N)}(\x, \vnu)$ is the fundamental Whittaker function,
which is normalized here as
$
\lim_{a \to 0} m^{(N)}(\x/a, a \vnu)
=e^{\vnu \cdot \x}
$
for $\x \in \W_N$.
Since $c_N(a \vnu) \simeq a^{N(N-1)/2} h_N(\vnu)$
as $a \to 0$, we have
\begin{equation}
\lim_{a \to 0} a^{N(N-1)/2} \psi^{(N)}_{a \vnu}(\x/a)
=\frac{\displaystyle{\det_{1 \leq j, k \leq N}[e^{\nu_j x_k}]}}{h_N(\vnu)} \quad
\mbox{for $\x, \vnu \in \W_N$}.
\label{eqn:tro3}
\end{equation}
Moreover, the density of Sklyanin measure (\ref{eqn:sN1})
has the asymptotics in $a \to 0$ as
$s_N(a \k) \simeq a^{N(N-1)} (h_N(\k))^2/\{(2 \pi)^N N! \}$,
and thus (\ref{eqn:QN1}) gives \cite{Kat12}
\begin{eqnarray}
\lim_{a \to 0} Q_N^{a}(t, \y|\x)
&=& \frac{1}{(2 \pi)^N N!}
\int_{\R^N} e^{-t|\k|^2/2}
\det_{1 \leq j, \ell \leq N} [e^{i k_j x_{\ell}}]
\det_{1 \leq j, \ell \leq N} [e^{-i k_j y_{\ell}}] d \k
\nonumber\\
&=& q_N(t, \y|\x), \quad
\x, \y \in \W_N, \quad t \in [0, \infty).
\label{eqn:tro4}
\end{eqnarray}
Then (\ref{eqn:pNnu3}) is concluded.
Note that the first equality in (\ref{eqn:tro4}) can be
interpreted in terms of the Slater determinants
used in quantum mechanics \cite{SMCR08}.

\vskip 0.5cm
%%%%%%%%%%%%%%%%%%%%%%%%%%%%%%%%%%%%%%%%%%%%%%%%%%%%%%
\noindent{\bf Acknowledgements} \quad
%%%%%%%%%%%%%%%%%%%%%%%%%%%%%%%%%%%%%%%%%%%%%%%%%%%%%%
%%%%%%%%%%%%%%%%%%%%%%%%%%%%%%%%%%%%%%%%%%%%%%%%%%
The present author would like to thank
T. Imamura and P. Graczyk for useful discussion
on diffusion processes with drifts.
A part of the present work was done
during the participation of the author 
in the EPSRC Symposium Workshop on
``Interacting particle systems, growth models
and random matrices"
at the university of Warwick (19-23 March 2012).
The author thanks N. O'Connell, J. Ortmann,
and J. Warren
for invitation to the workshop.
This work is supported in part by
the Grant-in-Aid for Scientific Research (C)
(No.21540397) of Japan Society for
the Promotion of Science.
%%%%%%%%%%%%%%%%%%%%%%%%%%%%%%%%%%%%%%%%%%%%%%%%%%

%%%%%%%%%%%%%%%%%%%%%%%%%%%%%%%%%%%%%%%%%%%%%%%%%%%%%%%%%%%%%
%%%%%%%%%%%%%%%Reference%%%%%%%%%%%%%%%%%%%%%%%%%%%%%%%%%%%%%
%%%%%%%%%%%%%%%%%%%%%%%%%%%%%%%%%%%%%%%%%%%%%%%%%%%%%%%%%%%%%

%%%%%%%%%%%%%%%%%%%%%%%%%%%%%%%%%%%%%%%%%%
\end{document}